\pdfoutput=1
\RequirePackage{ifpdf}
\ifpdf 
\documentclass[pdftex]{sigma}
\else
\documentclass{sigma}
\fi

\usepackage{tikz}

\begin{document}

\allowdisplaybreaks

\newcommand{\arXivNumber}{2007.13277}

\renewcommand{\PaperNumber}{134}

\FirstPageHeading

\ShortArticleName{Knot Complement, ADO Invariants and their Deformations for Torus Knots}

\ArticleName{Knot Complement, ADO Invariants\\ and their Deformations for Torus Knots}

\Author{John CHAE}

\AuthorNameForHeading{J.~Chae}

\Address{Univeristy of California Davis, Davis, USA}
\Email{\href{mailto:yjchae@ucdavis.edu}{yjchae@ucdavis.edu}}

\ArticleDates{Received August 20, 2020, in final form December 09, 2020; Published online December 15, 2020}

\Abstract{A relation between the two-variable series knot invariant and the Akutsu--Degu\-chi--Ohtsuki (ADO) invariant was conjectured recently. We reinforce the conjecture by presenting explicit formulas and/or an algorithm for particular ADO invariants of torus knots obtained from the series invariant of complement of a knot. Furthermore, one parameter deformation of ${\rm ADO}_3$ polynomial of torus knots is provided.}


\Classification{57K14; 57K16; 81R50}

\section{Introduction}\label{section1}

Categorification of link invariants has been a source of fruitful interactions between physics and low dimensional topology over the past decades (see \cite{Gukov,NO,Webster} for reviews). Since the advent of the Khovanov homology~\cite{Kh}, which categorifies the Jones polynomials of links, there has been constructions of other homological theories, for example, knot Floer homology~\cite{OS,JR}, Khovanov--Rozansky homology~\cite{KR} and HOMFLY homology~\cite{KR2} that categorify the well-known link polynomials: Alexander, $\mathfrak{sl}(N)$-invariants and HOMFLY polynomial, respectively. Not only has the categorification deepened the conceptual aspects of links, but it has also provided a more powerful machinery to compute higher structural invariants beyond polynomial invariants. Furthermore, these advancements have inspired new directions in physics, which resulted in physical realizations of the link homologies. Beginning from knot Floer homology, its physical interpretation was found in \cite{DGR}. A physical realization of Khovanov homology and Khovanov--Rozansky homology was first provided using topological string theory in \cite{GSV}; additionally, through the conifold transition, existence of the HOMFLY homology was predicted as well. In the case of Khovanov homology, a different physical system involving D-branes was achieved in~\cite{W2}. For Kauffman homology, its physical construction exemplified the role of orientifolds~\cite{GW}. Even knot homology based on an exceptional Lie algebra admits a physical description~\cite{EG} (see Table~\ref{table1} for summary).
\begin{table}[h!]\centering\small
\renewcommand{\arraystretch}{1.1}
\begin{tabular}{ |p{3.5cm}||p{4.5cm}|p{5cm}|} \hline
 \textbf{Polynomial} & \textbf{Homology} & \textbf{Physical realization}\\
 \hline
 Alexander & $\mathfrak{sl}(1\vert 1)$ knot Floer homology & M5-M2 branes on the deformed\newline conifold \\
Jones & $\mathfrak{sl}(2)$ Khovanov homology & M5-M2 branes on the deformed\newline conifold or D3-NS5 brane system \\
$\mathfrak{sl}(N)$-invariants & $\mathfrak{sl}(N)$ Khovanov--Rozansky\newline homology & M5-M2 branes on the deformed\newline conifold \\
HOMFLY & HOMFLY homology & M5-M2 branes on the resolved\newline conifold \\		
$\mathfrak{so}(n)/\mathfrak{sp}(n)$-invariants\newline \& Kauffman & Kauffman homology & D4-brane \& orientifold system\newline on the resolved conifold \\
hyperpolynomial & $e_6$ homology & M5-M2 branes on the resolved\newline conifold \\
\hline				
\end{tabular}
\caption{A summary of link invariants and their physical realizations. Choice of an orientifold type determines $\mathfrak{so}(n)$ or $\mathfrak{sp}(n)$ Lie algebra. Applications of $S$- and $T$-dualities are necessary to the latter brane system in the case of Khovanov homology (for details see~\cite{W2}).}\label{table1}
\end{table}

In recent years, a physical approach to categorification of the Witten--Reshetikhin--Tureav (WRT)-invariant of 3-manifold~\cite{RT2,RT,W1}, namely, homological blocks $\hat{Z}(q)$~\cite{GPPV,GPV} inspired a~new kind of invariant for a complement of a knot~\cite{GM}. This knot invariant denoted as $F_{K}$ is a~two-variable series that emerges from $\hat{Z}(q)$:
\begin{gather*}
F_{K}(x,q) : = \hat{Z}_{0} \big(M^{3}_{K}; x^{1/2},n,q\big),\qquad |q| < 1,
\end{gather*}
where $M^{3}_{K}$ is a complement of a knot $K$ in a closed oriented 3-manifold $M^3$, $n \in \mathbb{Z}$, $c \in \mathbb{Z}_{+}$ and~$\Delta \in \mathbb{Q}$.\footnote{The r.h.s.\ of the definition of $F_{K}$ is a two-variable version of $\hat{Z}_{b}(q)$.} Physical interpretation of $F_{K}$ is that it counts BPS states of a 3d $\mathcal{N}=2$ supersymmetric theory $T\big[M^{3}_{K}\big]$ on the knot complement, which arises from the \textit{integrality} of the coefficients of $F_{K}$-series. This in turn originates from the appearance of dimension of BPS Hilbert space of~$T[Y]$ in the q-series $\hat{Z}(Y,q)$ for a generic 3-manifold $Y$. Furthermore, this Hilbert space is identified with a conjectured triply graded three manifold homology $\mathcal{H}^{i,j}_{\rm BPS}(Y;b)$ whose (graded) Euler characteristic is
\begin{gather*}
\hat{Z}_{b}[Y;q]= \sum_{i,j} (-1)^i q^j \dim \mathcal{H}^{i,j}_{\rm BPS}(Y;b) \in 2^{-c} q^{\Delta} \mathbb{Z} [[q]],\qquad |q| < 1.
\end{gather*}
The WRT-invariant of $Y$ is recovered from $\hat{Z}_{b}[Y;q]$ as q goes to a root of unity (for details see \cite[Section~2]{GPPV}).

Among mathematical developments of $F_{K}$~\cite{S, P2,P1}, evidence for a relationship between~$F_{K}$ and the ADO link invariant~\cite{ADO} have been discovered in~\cite{GHNPPS}. This relation is conjectured to hold for all knots and for any roots of unity:

\begin{conjecture}[{\cite[Conjecture 3]{GHNPPS}}]\label{conjecture1} For any knot $K$ in $S^3$,
\begin{gather*}
F_K (x, q)|_{q=\zeta_{p}} = \big( x^{1/2}-x^{-1/2} \big) \frac{ {\rm ADO}_{p}(K; x,\zeta_{p})}{\Delta_{K}(x^p)},\qquad \zeta_{p}={\rm e}^{{\rm i} 2 \pi/p},\qquad p \in \mathbb{Z}_{+}.
\end{gather*}
This conjecture was verified for specific values of p for the right-handed trefoil and the figure eight knots~{\rm \cite{GHNPPS}}. Another advancement was an introduction of a refinement of $F_K(x,q)$~{\rm \cite{EGGKPS}}. It was shown that $F_K(x,q)$ admits two parameter deformations through the superpolynomial~{\rm \cite{DGR,FGSS}}. This led to a generalization of the above conjecture.
\end{conjecture}

\begin{conjecture}[{\cite[Conjecture 4]{EGGKPS}}]\label{conjecture2} For any knot $K$ in $S^3$, there exists a $t$-deformation of the symmetric ${\rm ADO}_{p}$-polynomial of $K$ for ${\rm SU}(N)$,
\begin{gather*}
{\rm ADO}_{p}^{{\rm SU}(N)}[K; x,t] : = \big( \Delta_K \big(x^p, -(-t)^p\big) \big)^{N-1} \lim_{q \rightarrow {\rm e}^{{\rm i}2\pi /p}}F_K(x,q,a=-q^N/t,t),\qquad p \in \mathbb{Z}_{+},
\end{gather*}
and $t=-1$ specialization reduces to the original ${\rm ADO}_p[K;x]$ $($up to rescaling of~$x)$.
\end{conjecture}

The rest of the paper is organized as follows. In Section~\ref{section2} we briefly review the series invariant for a knot complement and the ADO invariants. In Section~\ref{section3} we present the explicit formulas and/or an algorithm for the ${\rm ADO}_{3}$ and ${\rm ADO}_{4}$ polynomials for a particular class of torus knots. Furthermore, one parameter deformation of ${\rm ADO}_3$ invariants for torus knots is discussed.

\section{Background}\label{section2}
\subsection{Two-variable series knot complement}\label{section2.1}

A series invariant $F_K$ for a complement of a knot $M^{3}_K$ was introduced in \cite{GM}. It has various properties such as the gluing formula and the (Dehn) surgery formula. This knot invariant $F_K$ takes the form\footnote{Implicitly, there is a choice of group; originally, the group used is ${\rm SU}(2)$.}
\begin{gather*}
F_K(x,q)= \frac{1}{2} \sum_{\substack{m \geq 1 \\ m \ \text{odd}}}^{\infty} \big(x^{m/2}-x^{-m/2}\big)f_{m}(q) \in \frac{1}{2^{c}} q^{\Delta} \mathbb{Z}\big[x^{\pm 1/2}\big]\big[\big[q^{\pm 1}\big]\big],
\end{gather*}
where $f_{m}(q)$ are Laurent series with integer coefficients, $c \in \mathbb{Z}_{+}$ and $\Delta \in \mathbb{Q}$. Moreover, $x$-variable is associated to the relative ${\rm Spin}^c \big(M^{3}_K, T^2\big)$-structures, which is affinely isomorphic to $H^2\big(M^{3}_K, T^2 ; \mathbb{Z}\big) \cong H_1\big(M^{3}_K;\mathbb{Z}\big)$; it has an infinite order, which is reflected as a series in $F_K$. For applications, some classes of knots have been analyzed~\cite{GM,P2}. One of them is a class of torus knots, which is relevant for our purpose. Hence we display $F_{K}$ for the right-handed torus knots $T(s,t)$, $s,t>1$ with $\operatorname{gcd}(s,t)=1$~\cite{GM}.
\begin{gather*}
F_{T(s,t)}(x,q)= \frac{1}{2}q^{(s-1)(t-1)/2} \sum_{\substack{m \geq 1 \\ m \ \text{odd}}}^{\infty} \epsilon(s,t)_{m} \big(x^{m/2}-x^{-m/2}\big) q^{\frac{m^2-(st-s-t)^2}{4st}},
\\
\epsilon(s,t)_{m}= \begin{cases}
 -1,& m \equiv st+s+t\quad \text{or}\quad st-s-t \quad (\text{mod} \ 2st),\\
 \phantom{+}1, & m \equiv st+s-t \quad\text{or}\quad st-s+t \quad (\text{mod} \ 2st), \\
	 \phantom{+}0, & \text{otherwise}.
\end{cases}
\end{gather*}

Prior to $F_K$'s potential relation to the (original) ADO invariant, it was proposed that $F_K$ possess similar characteristics of $\mathfrak{sl}(2)$-colored Jones polynomial through the Melvin--Morton--Rozansky conjecture~\cite{MM,R} (proven in~\cite{BG}), and the quantum volume conjecture~\cite{G,Gukov2}:

\begin{conjecture}[{\cite[Conjecture 1.5]{GM}}]\label{conjecture3} For a knot $K \subset$ $S^3$, the asymptotic expansion of the knot invariant $F_{K}\big(x,q={\rm e}^{\hbar}\big)$ about $\hbar =0$ coincides with the Melvin--Morton--Rozansky expansion of the colored Jones polynomial in the large color limit:
\begin{gather*}
\frac{F_{K}\big(x,q={\rm e}^{\hbar}\big)}{x^{1/2}-x^{-1/2}} = \sum_{r=0}^{\infty} \frac{P_{r}(x)}{\Delta_K(x)^{2r+1}}\hbar^r,
\end{gather*}
where $x=q^{n\hbar}$ is fixed, n is the color of $K$, $P_{r}(x) \in \mathbb{Q} \big[x^{\pm 1}\big]$, $P_{0}(x)=1$ and $\Delta_K (x)$ is the Alexander polynomial of~$K$.
\end{conjecture}

\begin{conjecture}[{\cite[Conjecture 1.6]{GM}}] For any knot $K \subset$ $S^3$, the normalized series $f_{K}(x,q)$ satisfies a linear recursion relation generated by the quantum A-polynomial of $K$:
\begin{gather*}
\hat{A}_{K}(\hat{x},\hat{y},q) f_{K}(x,q) = 0,
\end{gather*}
where $f_{K}:=F_{K}(x,q)/\big(x^{1/2}-x^{-1/2}\big)$.
\end{conjecture}

\subsection{The ADO invariants of knots}\label{section2.2}
Colored generalization of the Alexander polynomial for framed colored and oriented knot (link) was introduced in~\cite{ADO}. This knot invariant(ADO invariant) is based on $(1,1)$-colored tangle diagram obtained by cutting the knot (or a component of a link). From this colored and oriented tangle diagram, the ADO invariant is constructed from a non-semisimple category of module over the unrolled quantum group $\mathcal{U}^{H}_{\zeta_{2r}}(\mathfrak{sl}_2(\mathbb{C}))$ together with the \textit{modified} quantum dimension ($r \in \mathbb{Z}_{\geq 2}$). We will employ the quantum algebra construction of the ADO invariants for verification of our results; the computational ingredients are summarized in Appendix~\ref{appendixB}. We give a concise review of the conceptual features of the construction~\cite{ADO,BDGG,SW}.

The first ingredient is the unrolled quantum group $\mathcal{U}^{H}_{\zeta_{2r}}({\mathfrak{sl}}_2(\mathbb{C}))$, which is a $\mathbb{C}$-algebra spe\-cialized at $q=\zeta_{2r}$; its generators and relations are
\begin{itemize}\itemsep=0pt
\item generators: $E$, $F$, $K$, $K^{-1}$, $H$,
\item relations:
\begin{gather*} KK^{-1}=K^{-1}K=1,\!\qquad KE=\zeta_{2r}^2 EK,\!\qquad KF=\zeta_{2r}^{-2} FK,\!\qquad [E,F]=\frac{K-K^{-1}}{\zeta_{2r}-\zeta_{2r}^{-1}},\\
KH=HK, \qquad [H,E]=2E, \qquad [H,F]=-2F, \qquad E^r=F^r=0.
\end{gather*}
\end{itemize}
This algebra possess a Hopf algebra structure:
\begin{gather*}
\Delta (E)= 1 \otimes E + E \otimes K, \qquad \epsilon(E)=0, \qquad \quad S(E)=-EK^{-1},
\\
\Delta (F)= K^{-1} \otimes F + F \otimes 1, \qquad \epsilon(F)=0, \qquad S(F)=-KF,
\\
\Delta (H)= 1 \otimes H + H \otimes 1, \qquad \epsilon(H)=0, \qquad S(H)=-H,
\\
\Delta (K)= K \otimes K, \qquad \epsilon(K)=1, \qquad S(K)=K^{-1},
\\
\Delta \big(K^{-1}\big)= K^{-1} \otimes K^{-1}, \qquad \epsilon\big(K^{-1}\big)=1, \qquad S\big(K^{-1}\big)=K.
\end{gather*}
The second element of the construction of the ADO invariant is a functor RT between a category of colored oriented tangle diagrams COD and a category Rep of representations of $\mathcal{U}^{H}_{\zeta_{2r}}({\mathfrak{sl}}_2(\mathbb{C}))$:
\begin{gather*}
{\rm RT} \colon \ {\rm COD} \longrightarrow {\rm Rep}.
\end{gather*}
The objects of COD are framed colored oriented $(1,1)$-tangle diagrams and morphisms are equivalence classes of the tangle diagrams whose equivalence relations are generated by the tangle moves (see \cite[Section~2]{ADO}). For the target category, its objects are vector spaces V and morphisms are linear maps between them. The image of the RT functor is ${\rm RT}(T)=\langle T\rangle \operatorname{Id}_{V} \in \operatorname{End}_{\mathbb{C}}(V)$, which enables to define
\begin{gather*}
{\rm ADO}(K)_{r} := d(V_{\alpha};r)\langle T\rangle,
\end{gather*}
where $V_{\alpha}$ is a vector space assigned to $K$ (or to an open component of a link\footnote{ADO invariant is independent of choice of a component of a link that is cut (for details see \cite[Section~5]{ADO}).}) and $d(V_{\alpha};r)$ is the modified quantum dimension,
\begin{gather*}
d(V_{\alpha};r) = -\zeta_{2r}^{\frac{1}{2}r(1-r)}\frac{\zeta_{2r}^{\alpha +1}-\zeta_{2r}^{-\alpha -1}}{\zeta_{2r}^{r\alpha}-\zeta_{2r}^{-r\alpha}},\qquad \alpha \in (\mathbb{C} \backslash \mathbb{Z})\cup ( r\mathbb{Z} -1).
\end{gather*}
This modified dimension replaces the usual quantum trace, which vanishes in this context. Moreover, it makes ${\rm ADO}(K)$ an isotopy invariant.

\section{The ADO invariants of torus knots}\label{section3}
Recently, evidence for a relation between $F_{K}$ at specific values of roots of unity and the ADO invariants were discovered for the (right-handed) trefoil, the figure eight and $5_2$ knots~\cite{GHNPPS}. Furthermore, this relation is conjectured to hold for any roots of unity and for all knots (Conjec\-ture~\ref{conjecture1}). Using the formula in Section~\ref{section2.1} and Conjecture~\ref{conjecture1}, close examination of torus knots $T(2,2s+1)$ at various values of s yields an explicit formula or an algorithm for ${\rm ADO}_3$ and ${\rm ADO}_4$ invariants of $T(2,2s+1)$, $s \in \mathbb{Z}_{+}$.

\subsection[The ADO3 invariants of T(2,2s+1)]{The $\boldsymbol{{\rm ADO}_3}$ invariants of $\boldsymbol{T(2,2s+1)}$}\label{section3.1}

The ${\rm ADO}_3$ invariants of $T(2,2s+1)$ are divided in three types depending on their coefficient pattern:
\begin{enumerate}\itemsep=0pt
	\item[1)] for $K=T(2,2s+1)= T(2,3), T(2,9),T(2,15),T(2,21),\dots$
		\begin{gather*}
	 {\rm ADO}_3(x) = \zeta_{3}x^{2s} + \zeta_{3}x^{2s-1} + \big(\zeta_{3}-\zeta_{3}^{-1}\big)x^{2s-2} - \zeta_{3}^{-1}x^{2s-3} -\zeta_{3}^{-1}x^{2s-4}\\	
\hphantom{{\rm ADO}_3(x) =}{} +\zeta_{3}x^{2s-6} + \zeta_{3}x^{2s-7} + \big(\zeta_{3}-\zeta_{3}^{-1}\big)x^{2s-8} - \zeta_{3}^{-1}x^{2s-9} -\zeta_{3}^{-1}x^{2s-10}+ \cdots\\
\hphantom{{\rm ADO}_3(x) =}{} + \big(\zeta_{3}-\zeta_{3}^{-1}\big) + (x \rightarrow 1/x),
	\end{gather*}
	
\item[2)] for $K=T(2,2s+1)= T(2,5), T(2,11),T(2,17),T(2,23),\dots$
		\begin{gather*}
	 {\rm ADO}_3(x) = \zeta_{3}^{-1}x^{2s} + \zeta_{3}^{-1} x^{2s-1} + \big(\zeta_{3}^{-1}-1\big) x^{2s-2} - x^{2s-3} - x^{2s-4} +\zeta_{3}^{-1}x^{2s-6} \\
\hphantom{{\rm ADO}_3(x) =}{} + \zeta_{3}^{-1} x^{2s-7}+ \big(\zeta_{3}^{-1}-1\big) x^{2s-8} - x^{2s-9} - x^{2s-10} + \cdots -1 + (x \rightarrow 1/x),
		\end{gather*}
	
\item[3)] for $K=T(2,2s+1)= T(2,7), T(2,13),T(2,19),T(2,25),\dots$
		\begin{gather*}
	 {\rm ADO}_3(x) = x^{2s} + x^{2s-1} + (1-\zeta_{3}) x^{2s-2} - \zeta_{3} x^{2s-3} - \zeta_{3} x^{2s-4}+ x^{2s-6} + x^{2s-7} \\
\hphantom{{\rm ADO}_3(x) =}{} + (1-\zeta_{3}) x^{2s-8} - \zeta_{3} x^{2s-9} - \zeta_{3} x^{2s-10}+ \cdots +1 + (x \rightarrow 1/x).
	\end{gather*}
\end{enumerate}
All the explicit x terms are polynomials and power of x decreases by two after one cycle of a~coefficient combination. We next move onto ${\rm ADO}_4$ invariants, whose explicit formula can be obtained algorithmically.

\subsection[The algorithm for ADO4 invariants of T(2,2s+1)]{The algorithm for $\boldsymbol{{\rm ADO}_4}$ invariants of $\boldsymbol{T(2,2s+1)}$}\label{section3.2}

Explicit formulas for ${\rm ADO}_{4}$ invariants of $T(2,2s+1)$ for $s \in \mathbb{Z}_{\geq 7}$ are constructed inductively. This subclass of torus knots are divided into four sets and each set has its own seed ${\rm ADO}_4[T(2,2s+1)]$ together with a pattern of coefficients that generates the invariant for higher values of $2s+1$. We present an algorithm for obtaining explicit expressions.

\textbf{The algorithm:}
\begin{enumerate}\itemsep=0pt
	\item Beginning with $x^{3s}$, write a polynomial with coefficients $c_{i}$ following one of the four patterns (shown below) that $T(2,2s+1)$ belong to
	\begin{gather*}
	c_{3s} x^{3s}\! + c_{3s-1} x^{3s-1}\! + c_{3s-2} x^{3s-2}\! + c_{3s-3} x^{3s-3}\! + c_{3s-4} x^{3s-4}\! + c_{3s-5} x^{3s-5},\!\qquad c_{n} \in \mathbb{C}.
	\end{gather*}
	
	\item Add a polynomial starting with $x^{3s-8}$ with exponent pattern and coefficients given by
	${\rm ADO}_4[T(2,2s-7)]$.
	
	\item Remaining polynomials are determined by a mirror reflection of coefficients across the last term in the previous step beginning from the second last term. Furthermore, adjust of the exponents of the variable~$x$ following the pattern of Step~2 until a constant term is reached.
	
	\item Use the Weyl symmetry to obtain $1/x$ terms.
	
\end{enumerate}
As a consequence of the normalization factor $\big(x^{1/2}-x^{-1/2}\big)$ in Conjecture~\ref{conjecture1}, we obtain the symmetric version of ADO invariants. Their coefficients ${c_{n}}$ are divided into four types:
\begin{enumerate}\itemsep=0pt
	\item[1)] $-{\rm i}$, $-{\rm i}$, $-{\rm i}-1$, $-{\rm i}-1$, $-1$, $-1$ \text{for} $\{ T(2,7), T(2,15), T(2,23),\dots \} $,
	
	\item[2)] $1$, $1$, $1-{\rm i}$, $1-{\rm i}$, $-{\rm i}$, $-{\rm i}$ \text{for} $\{ T(2,9), T(2,17), T(2,25),\dots \} $,
	
	\item[3)] ${\rm i}$, ${\rm i}$, ${\rm i}+1$, ${\rm i}+1$, $1$, $1$ \text{for} $\{ T(2,11), T(2,19),\dots \} $,
	
	\item[4)] $-1$, $-1$, $-1+{\rm i}$, $-1+{\rm i}$, ${\rm i}$, ${\rm i}$ \text{for} $\{ T(2,13), T(2,21),\dots \} $,
\end{enumerate}
where the semicolon means that the next term has a power of $x$ lowered by three. The coefficients of the first and the third sets differ by signs as well as the second and the fourth sets. ${\rm ADO}_4$ polynomial of the first knot in each set is a seed for the next knot in the set. This pattern continues for all the subsequent knots in each set. The fundamental seed invariants can be easily computed using the torus knot formula $F_{T(s,t)}$ in Section~\ref{section2.1}
\begin{gather*}
	 {\rm ADO}_4[T(2,7)] = -{\rm i} x^{9} - {\rm i} x^{8} + (-1-{\rm i})x^{7} + (-1-{\rm i})x^{6} - x^{5} - x^{4} - {\rm i}x^{2} - {\rm i}2 x\\
\hphantom{{\rm ADO}_4[T(2,7)] =}{} + 1-{\rm i}2 + (x \rightarrow 1/x ),\\
{\rm ADO}_4[T(2,9)] = x^{12} + x^{11} + (1-{\rm i})x^{10} + (1-{\rm i})x^{9} - {\rm i} x^8 - {\rm i} x^7 + x^4 -{\rm i} x^2 -{\rm i}2 x\\
\hphantom{{\rm ADO}_4[T(2,9)] =}{} - 1 -{\rm i}2+ (x \rightarrow 1/x ),	
\\
{\rm ADO}_4[T(2,11)]= {\rm i}x^{15} + {\rm i}x^{14} + (1+{\rm i})x^{13} + (1+{\rm i})x^{12} + x^{11} + x^{10} + {\rm i}x^7 + {\rm i}x^6\\
\hphantom{{\rm ADO}_4[T(2,11)]=}{}	+ (1+{\rm i})x^5 + (1+{\rm i}2)x^4 + (1+{\rm i})x^3 + {\rm i}x^2 + ({\rm i}-1)x - 1 + (x \rightarrow 1/x ),
\\
{\rm ADO}_4[T(2,13)]= -x^{18} - x^{17} + (-1+{\rm i})x^{16} + (-1+{\rm i})x^{15} + {\rm i} x^{14} + {\rm i} x^{13} -x^{10} - x^{9}\\
\hphantom{{\rm ADO}_4[T(2,13)]=}{} + (-1+{\rm i})x^{8} + (-1+{\rm i})x^{7} + {\rm i} x^{6} + (1+{\rm i}) x^{5} + x^4 + (1+{\rm i})x^3+ {\rm i}x^2\\
\hphantom{{\rm ADO}_4[T(2,13)]=}{}	+ ({\rm i}-1)x-1+{\rm i}2 + (x \rightarrow 1/x ).
\end{gather*}
For completeness, we display the ${\rm ADO}_4$ polynomials of $T(2,3)$~\cite{GHNPPS} and $T(2,5)$
\begin{gather*}
{\rm ADO}_4[T(2,3)]= {\rm i}x^3 + {\rm i}x^2 + (1+{\rm i})x + 1+{\rm i}2 + (x \rightarrow 1/x ),
\\
{\rm ADO}_4[T(2,5)]= -x^6 -x^5 + (-1+{\rm i})x^4 + (-1+{\rm i})x^3 + {\rm i} x^2 + (1+{\rm i})x + 1 + (x \rightarrow 1/x ).
\end{gather*}

\subsection{Examples}\label{section3.3}
Let us demonstrate the algorithm through examples. For $T(2,15)$ in the first set, the first step of the algorithm yields
\begin{gather*} \text{Step 1}= -{\rm i} x^{21} -{\rm i}x^{20}+ (-1-{\rm i})x^{19} + (-1-{\rm i})x^{18}-x^{17}-x^{16}. \end{gather*}
Next step is to use the coefficients from the seed ${\rm ADO}_4[T(2,7)]$ but its powers of x are adjusted appropriately
\begin{gather*}
 \text{Step 2}= -{\rm i} x^{21} -{\rm i}x^{20}+ (-1-{\rm i})x^{19} + (-1-{\rm i})x^{18}-x^{17}-x^{16}\\
\hphantom{\text{Step 2}=}{} -{\rm i} x^{13} - {\rm i} x^{12} + (-1-{\rm i})x^{11} + (-1-{\rm i})x^{10} - x^{9} - x^{8} - {\rm i}x^{6} - {\rm i}2 x^{5} + (1-{\rm i}2)x^{4}.
\end{gather*}
Since the above expression ends in $(1-{\rm i}2)x^{4}$, we need to reflect the coefficients about this term until a constant term is reached. This results in
\begin{gather*}
\text{Step 3}= -{\rm i} x^{21} -{\rm i}x^{20}+ (-1-{\rm i})x^{19} + (-1-{\rm i})x^{18}-x^{17}-x^{16} -{\rm i} x^{13} - {\rm i} x^{12}+ (-1-{\rm i})x^{11}\\
\hphantom{\text{Step 3}=}{} + (-1-{\rm i})x^{10} - x^{9} - x^{8} - {\rm i}x^{6} - {\rm i}2 x^{5} + (1-{\rm i}2)x^{4}- {\rm i}2 x^{3} - {\rm i}x^{2} - 1.
\end{gather*}
The application of the last step leads to
\begin{gather*}
{\rm ADO}_4[T(2,15)]= -{\rm i} x^{21} -{\rm i}x^{20}+ (-1-{\rm i})x^{19} + (-1-{\rm i})x^{18}-x^{17}-x^{16} -{\rm i} x^{13} - {\rm i} x^{12}\\
\hphantom{{\rm ADO}_4[T(2,15)]=}{} + (-1-{\rm i})x^{11} + (-1-{\rm i})x^{10} - x^{9} - x^{8} - {\rm i}x^{6} - {\rm i}2 x^{5} + (1-{\rm i}2)x^{4}\\
\hphantom{{\rm ADO}_4[T(2,15)]=}{} - {\rm i}2 x^{3} - {\rm i}x^{2} - 1 + (x \rightarrow 1/x ).
\end{gather*}
As a consistency check, $F_{T(2,15)}(x,q=\zeta_4)$ obtained from the ${\rm ADO}_4[T(2,15)]$ using Conjecture~\ref{conjecture1} agrees with the direct computation of $F_{T(2,15)}(x,q=\zeta_4)$ from Section~\ref{section2.1}.

For $T(2,17)$ in the second set, the seed invariant is ${\rm ADO}_4[T(2,9)]$ and application of the first and second steps produce
\begin{gather*}
\text{Step 1}= x^{24} + x^{23}+ (1-{\rm i})x^{22} + (1-{\rm i})x^{21}- {\rm i} x^{20}- {\rm i} x^{19},
\\
 \text{Step 2}= x^{24} + x^{23}+ (1-{\rm i})x^{22} + (1-{\rm i})x^{21}- {\rm i} x^{20}- {\rm i} x^{19}+ x^{16} + x^{15}\\
\hphantom{\text{Step 2}=}{}+ (1-{\rm i})x^{14} + (1-{\rm i})x^{13}- {\rm i} x^{12}- {\rm i} x^{11} + x^{8} - {\rm i} x^{6} -{\rm i}2 x^{5} + (-1-{\rm i}2 )x^{4}.
\end{gather*}
After the reflection about $x^{4}$-term
\begin{gather*}
\text{Step 3}= x^{24} + x^{23}+ (1-{\rm i})x^{22} + (1-{\rm i})x^{21}- {\rm i} x^{20}- {\rm i} x^{19} + x^{16} + x^{15}+ (1-{\rm i})x^{14} \\
\hphantom{\text{Step 3}=}{} + (1-{\rm i})x^{13}- {\rm i} x^{12}- {\rm i} x^{11} + x^{8} - {\rm i} x^{6} -{\rm i}2 x^{5} + (-1-{\rm i}2 )x^{4} - {\rm i}2 x^{3} - {\rm i} x^{2} + 1.
\end{gather*}
The last step results in
\begin{gather*}
{\rm ADO}_4[T(2,17)]= x^{24} + x^{23}+ (1-{\rm i})x^{22} + (1-{\rm i})x^{21}- {\rm i} x^{20}- {\rm i} x^{19}+ x^{16} + x^{15} \\
\hphantom{{\rm ADO}_4[T(2,17)]=}{} +(1-{\rm i})x^{14}+ (1-{\rm i})x^{13}- {\rm i} x^{12}- {\rm i} x^{11} + x^{8} - {\rm i} x^{6} -{\rm i}2 x^{5} + (-1-{\rm i}2 )x^{4}\\
\hphantom{{\rm ADO}_4[T(2,17)]=}{} - {\rm i}2 x^{3} - {\rm i} x^{2} + 1+ (x \rightarrow 1/x ).
\end{gather*}
One can verify that $F_{T(2,17)}(x,q=\zeta_4)$ obtained from ${\rm ADO}_4[T(2,17)]$ matches with the result (at $q=\zeta_4$) of the direct method from Section~\ref{section2.1}.

In the third set, the seed for $T(2,19)$ is ${\rm ADO}_4[T(2,11)]$. Applying the first two steps yields
\begin{gather*}
 {\rm i} x^{27}+ {\rm i} x^{26} + (1+{\rm i})x^{25} + (1+{\rm i})x^{24} + x^{23} + x^{22} + {\rm i} x^{19}+ {\rm i} x^{18} + (1+{\rm i})x^{17} + (1+{\rm i})x^{16}\\
\quad{}+ x^{15} + x^{14} + {\rm i} x^{11}+ {\rm i} x^{10} + (1+{\rm i})x^{9} + (1+{\rm i}2)x^{8} + (1+{\rm i}) x^{7} + {\rm i} x^{6} + (-1+{\rm i})x^{5} - x^{4}.
\end{gather*}
The last two steps produce
\begin{gather*}
{\rm ADO}_4[T(2,19)]= {\rm i} x^{27}+ {\rm i} x^{26} + (1+{\rm i})x^{25} + (1+{\rm i})x^{24} + x^{23} + x^{22} + {\rm i} x^{19}+ {\rm i} x^{18}\\
\hphantom{{\rm ADO}_4[T(2,19)]=}{} + (1+{\rm i})x^{17} + (1+{\rm i})x^{16} + x^{15} + x^{14} + {\rm i} x^{11}+ {\rm i} x^{10} + (1+{\rm i})x^{9}\\
\hphantom{{\rm ADO}_4[T(2,19)]=}{} + (1+{\rm i}2)x^{8} + (1+{\rm i}) x^{7} + {\rm i} x^{6} + (-1+{\rm i})x^{5} - x^{4} + (-1+{\rm i})x^{3} + {\rm i} x^{2} \\
\hphantom{{\rm ADO}_4[T(2,19)]=}{} + (1+{\rm i}) x + 1+{\rm i}2 + (x \rightarrow 1/x ).
\end{gather*}
Similarly, ${\rm ADO}_4[T(2,21)]$ can be obtained using ${\rm ADO}_4[T(2,13)]$
\begin{gather*}
{\rm ADO}_4[T(2,21)]= - x^{30} - x^{29}+ (-1+{\rm i})x^{28} + (-1+{\rm i})x^{27} + {\rm i} x^{26} + {\rm i} x^{25} - x^{22} - x^{21} \\
 \hphantom{{\rm ADO}_4[T(2,21)]=}{}+ (-1+{\rm i})x^{20} + (-1+{\rm i})x^{19} + {\rm i} x^{18}+ {\rm i} x^{17} - x^{14} - x^{13} + (-1+{\rm i}) x^{12}\\
\hphantom{{\rm ADO}_4[T(2,21)]=}{} + (-1+{\rm i})x^{11} + {\rm i} x^{10} + (1+{\rm i}) x^{9} + x^{8} + (1+{\rm i})x^{7} + {\rm i}x^{6} + (-1+{\rm i})x^{5}\\
\hphantom{{\rm ADO}_4[T(2,21)]=}{} + (-1+{\rm i}2)x^{4} + (-1+{\rm i})x^{3} + {\rm i} x^{2} + (1+{\rm i})x + 1 + (x \rightarrow 1/x ).
\end{gather*}
Formulas for ${\rm ADO}_4$ invariants become lengthy as the winding number along the longitude of a~torus increases so their expressions are recorded in Appendix~\ref{appendixA}. We move onto the deformation of the ADO polynomial.

\subsection[Deformed ADO3 invariants of T(2,2s+1)]{Deformed $\boldsymbol{{\rm ADO}_3}$ invariants of $\boldsymbol{T(2,2s+1)}$}\label{section3.4}

A link between superpolynomial defined in~\cite{DGR} and $F_K$ was discovered in~\cite{EGGKPS}. Specifically, two parameter refinement $F_K(x,q,a,t)$ was introduced, which motivated to define $t$-deformed ADO polynomial. This deformation introduces one more variable to the original ADO polynomial~${\rm ADO}(x,t)$; as a~consequence, it is a colored version of the $t$-deformed Alexander polynomial~$\Delta(x,t)$ that can distinguish chirality of torus knots. In this Subsection we present $t$-deformed version of~${\rm ADO}_3$ polynomials for $T(2,2s+1)$ knots.

Reduced superpolynomial for the right-handed torus knots carrying symmetric representation~$S^r$ of ${\rm SU}(N)$ is stated in~\cite{FGSS}:
\begin{gather*}
\mathcal{P}_{S^r}[T(2,-(2s+1));q,a,t]= \left( \frac{a}{q} \right)^{pr} \sum_{k_{1}=0}^{r}\sum_{k_{2}=0}^{k_1}\cdots \sum_{k_{s}=0}^{k_{s-1}} q^{(2r+1)(k_1 + \cdots + k_s) - \sum\limits_{i=1}^{s} k_{i-1}k_{i}} t^{2(k_1 + \cdots + k_s)}\\
 \hphantom{\mathcal{P}_{S^r}[T(2,-(2s+1));q,a,t]=}{} \times \frac{\big(q^r;q^{-1}\big)_{k_1} (-at/q;q)_{k_1}}{(q;q)_{k_1}} \left[ \begin{matrix} k_1 \\ k_2 \end{matrix}\right]_{q} \cdots \left[ \begin{matrix} k_{s-1} \\ k_s \end{matrix}\right]_{q},
\\
(w;q)_m := \prod_{i=1}^{m} \big( 1- wq^{i-1} \big), \qquad \left[ \begin{matrix} w \\ n \end{matrix}\right]_q := \frac{(q;q)_w}{(q;q)_{n}(q;q)_{w-n} },
\end{gather*}
where $s \in \mathbb{Z}_{+}$, $r$ is the dimension of $S^r$ and $k_{0} \equiv r$. Note that the convention for the left-handed torus knot in \cite{EGGKPS} is $T(2,2s+1)$ for $s \in \mathbb{Z}_{+}$, which is opposite of the convention used in this article. In \cite{EGGKPS}, it was shown that $\mathcal{P}_{S^r}$ can be converted into a two parameter deformation of $F_K$ by replacing $q^r$ by x and dropping the overall factor $(a/q)^{pr}$:
\begin{gather*}
F_{T(2,-(2s+1))}(x,q,a,t)= \sum_{k_{1}=0}^{\infty}\sum_{k_{2}=0}^{k_1}\cdots \sum_{k_{s}=0}^{k_{s-1}} x^{2(k_1 + \cdots + k_s)-k_1} q^{(k_1 + \cdots + k_s) - \sum\limits_{i=2}^{s} k_{i-1}k_{i}} t^{2(k_1 + \cdots + k_s)}\\
\hphantom{F_{T(2,-(2s+1))}(x,q,a,t)=}{} \times \frac{\big(x;q^{-1}\big)_{k_1} (-at/q;q)_{k_1}}{(q;q)_{k_1}} \left[ \begin{matrix} k_1 \\ k_2 \end{matrix}\right]_{q} \cdots \left[ \begin{matrix} k_{s-1} \\ k_s \end{matrix}\right]_{q}.
\end{gather*}
Fixing $a=q^N$ and $t=-1$, $F_{K}(x,q,a,t)$ becomes the original $F_K(x,q)$ for torus knots.\footnote{Specifically, additional manipulations are needed to arrive at $F_K(x,q)$ for torus knots \cite[Section~5.2]{EGGKPS}.} Different specialization of a, namely, $a=-t^{-1}$ yields a refined Alexander polynomial~\cite{EGGKPS},
\begin{gather*}
F_K\big(x,q,-t^{-1},t\big)= \Delta_K(x,t).
\end{gather*}
Using Conjecture~\ref{conjecture2}, a refined ${\rm ADO}_3$ polynomial for $T(2,2s+1)$, $s \in \mathbb{Z}_{+}$ is
\begin{gather*}
{\rm ADO}_{3}[T(2,2s+1);x,t]= (tx)^{2s} + \frac{\zeta_{3}^{-1}}{t}(tx)^{2s-1} + \left( \frac{\zeta_{3}}{t^2} - \zeta_{3}^{-1} \right) (tx)^{2s-2} - \frac{\zeta_{3}}{t}(tx)^{2s-3}
\\
\hphantom{{\rm ADO}_{3}[T(2,2s+1);x,t]=}{} -\frac{1}{t^2}(tx)^{2s-4}+ (tx)^{2s-6} + \frac{\zeta_{3}^{-1}}{t}(tx)^{2s-7} + \left( \frac{\zeta_{3}}{t^2} - \zeta_{3}^{-1} \right) (tx)^{2s-8}\\
\hphantom{{\rm ADO}_{3}[T(2,2s+1);x,t]=}{} - \frac{\zeta_{3}}{t}(tx)^{2s-9} -\frac{1}{t^2}(tx)^{2s-10} + \cdots + O\left( \frac{1}{tx} \right),
\end{gather*}
where $O(1/tx)$-terms are determined by the $t$-deformed Weyl symmetry of the ${\rm ADO}_{p}$ invariant,
\begin{gather*}
{\rm ADO}^{{\rm SU}(2)}_{p} (1/x ,t ) = {\rm ADO}^{{\rm SU}(2)}_{p} \big( \zeta_{p}^{-2}t^{-2}x ,t \big).
\end{gather*}
The suppressed polynomial terms follow the same power and coefficient patterns of the previous terms. The three formulas for the original ${\rm ADO}_3[T(2,2s+1);x]$ coalesce into one formula through the $t$-deformation. We next present a few examples.

{\bf $\boldsymbol{K=T(2,5)}$.}
We start from $F_K(x,q,a,t)$ for $T(2,-5)$,
\begin{gather*}
F_{T(2,-5)}(x,q,a,t)= \sum_{k_{1}=0}^{\infty} \sum_{k_{2}=0}^{k_{1}} x^{2(k_{1}+k_{2})-k_{1}} q^{k_{1}+k_{2}-k_{1}k_{2}} t^{2(k_{1}+k_{2})} \frac{\big(x;q^{-1}\big)_{k_{1}} \big({-}\frac{a t}{q};q\big)_{k_{1}}}{(q;q)_{k_{1}}} \left[ \begin{matrix} k_1 \\ k_2 \end{matrix}\right]_{q}.
\end{gather*}
We next apply the mirror map to reverse the orientation of $K$,
\begin{gather*}
x \mapsto 1/x, \qquad q \mapsto 1/q, \qquad a \mapsto 1/a, \qquad t \mapsto 1/t.
\end{gather*}
Setting $a=-1/t$, we get a refined Alexander polynomial of $K$ (upon multiplication by an overall monomial),
\begin{gather*}
\Delta_K(x,t)= t^2 x^2+\frac{1}{t^2 x^2}-\frac{1}{t^2 x}-x+1.
\end{gather*}
Further fixing $t=-1$, it reduces to the Alexander polynomial of~$K$. Moreover, this refined polynomial possess the $t$-deformed Weyl symmetry for the refined Alexander polynomial,
\begin{gather*}
\Delta_K(1/x,t)=\Delta_K\big(x/t^2,t\big).
\end{gather*}
A refined ${\rm ADO}_3$ polynomial of $K$ is computed via Conjecture~\ref{conjecture2} as
\begin{gather*}
{\rm ADO}_{3}[T(2,5);x,t]= (tx)^4 + \frac{\zeta_{3}^{-1}}{t}(tx)^3 + \left( \frac{\zeta_{3}}{t^2} - \zeta_{3}^{-1} \right)(tx)^2 - \frac{\zeta_{3}}{t}(tx) -\frac{1}{t^2} - \frac{\zeta_{3}^{-1}}{t}\frac{1}{(tx)}
\\
\hphantom{{\rm ADO}_{3}[T(2,5);x,t]=}{}+ \left( \frac{1}{t^2} - \zeta_{3} \right) \frac{1}{(tx)^2} + \frac{\zeta_{3}^{-1}}{t}\frac{1}{(tx)^3} + \zeta_{3}\frac{1}{(tx)^4}.
\end{gather*}
This formula carries the $t$-deformed Weyl symmetry of the ${\rm ADO}_3$ invariant. Moreover, fixing $t=-1$ and rescaling $x \mapsto \zeta_{3}^2 x$, the refined polynomial becomes the original ${\rm ADO}_3$ polynomial,
\begin{gather*}
\zeta_{3}^{-1}x^4 + \zeta_{3}^{-1}x^3 + \big(\zeta_{3}^{-1}-1\big)x^2 - x -1 + (x \rightarrow 1/x ).
\end{gather*}

{\bf $\boldsymbol{K=T(2,7)}$.}
Two parameter deformation of $F_K$ for $T(2,-7)$ is
\begin{gather*}
F_{T(2,-7)}(x,q,a,t)= \sum_{k_{1}=0}^{\infty} \sum_{k_{2}=0}^{k_{1}} \sum_{k_{3}=0}^{k_{2}} x^{2 (k_{1}+k_{2}+k_{3})-k_{1}} q^{k_{1}+k_{2}+k_{3}-k_{1}k_{2}-k_{2} k_{3}} t^{2(k_{1}+k_{2}+k_{3})}\\
\hphantom{F_{T(2,-7)}(x,q,a,t)=}{} \times \frac{ \big(x;q^{-1}\big)_{k_{1}} \big({-}\frac{at}{q};q\big)_{k_{1}}}{(q;q)_{k_{1}}} \left[ \begin{matrix} k_1 \\ k_2 \end{matrix}\right]_{q} \left[ \begin{matrix} k_2 \\ k_3 \end{matrix}\right]_{q}.
\end{gather*}
A refined Alexander polynomial of $K$ having the refined Weyl symmetry is
\begin{gather*}
\Delta_{T(2,7)}(x,t)= -t^3 x^3-\frac{1}{t^3 x^3}+\frac{1}{t^3 x^2}+t x^2-t x-\frac{1}{t x}+\frac{1}{t}.
\end{gather*}
A refined ${\rm ADO}_3$ polynomial of $K$ is
\begin{gather*}
{\rm ADO}_{3}[T(2,7);x,t]= (tx)^6 + \frac{\zeta_{3}^{-1}}{t}(tx)^5 + \left( \frac{\zeta_{3}}{t^2} - \zeta_{3}^{-1} \right)(tx)^4 - \frac{\zeta_{3}}{t}(tx)^3 -\frac{1}{t^2}(tx)^2 +1
\\
\hphantom{{\rm ADO}_{3}[T(2,7);x,t]=}{} - \frac{\zeta^{-1}}{t^2}\frac{1}{(tx)^2} - \frac{\zeta}{t}\frac{1}{(tx)^3} + \left( \frac{\zeta^{-1}}{t^2}-1 \right) \frac{1}{(tx)^4} + \frac{\zeta}{t}\frac{1}{(tx)^5} +\frac{1}{(tx)^6}.
\end{gather*}
This polynomial possess the $t$-deformed Weyl symmetry of the ${\rm ADO}_3$ invariant and after specializing $t=-1$ and rescaling $x \mapsto \zeta_{3}^2 x$, it becomes{\samepage
\begin{gather*} x^6 + x^5 + (1-\zeta_{3} )x^4 -\zeta_{3} x^3 - \zeta_{3} x^2 + 1 + (x \rightarrow 1/x ), \end{gather*}
which is the original ${\rm ADO}_3$ polynomial for~$K$.}

{\bf $\boldsymbol{K=T(2,9)}$.} A refined Alexander polynomial of $K$ carrying the refined Weyl symmetry is
\begin{gather*}
\Delta_{T(2,9)}(x,t)= t^4 x^4+\frac{1}{t^4 x^4}-\frac{1}{t^4 x^3}-t^2 x^3+t^2 x^2+\frac{1}{t^2 x^2}-\frac{1}{t^2 x}-x+1.
\end{gather*}
A refined ${\rm ADO}_3$ polynomial of $K$ is
\begin{gather*}
{\rm ADO}_{3}[T(2,9);x,t]= (tx)^8 + \frac{\zeta_{3}^{-1}}{t}(tx)^7 + \left( \frac{\zeta_{3}}{t^2} - \zeta_{3}^{-1} \right)(tx)^6 - \frac{\zeta_{3}}{t}(tx)^5 -\frac{1}{t^2}(tx)^4
\\
\hphantom{{\rm ADO}_{3}[T(2,9);x,t]=}{} +(tx)^2 + \frac{\zeta_{3}^{-1}}{t}(tx) + \left( \frac{\zeta_{3}}{t^2} - \zeta_{3}^{-1} \right) + \frac{1}{t}\frac{1}{tx} + \zeta_{3}^{2}\frac{1}{(tx)^2} - \frac{\zeta_{3}}{t^2}\frac{1}{(tx)^4}
\\
\hphantom{{\rm ADO}_{3}[T(2,9);x,t]=}{}- \frac{1}{t}\frac{1}{(tx)^5} + \left( \frac{\zeta_{3}}{t^2} - \zeta_{3}^{-1} \right) \frac{1}{(tx)^6}
 + \frac{1}{t}\frac{1}{(tx)^7} + \zeta_{3}^{2}\frac{1}{(tx)^8}.
\end{gather*}
This polynomial is invariant under the refined Weyl symmetry of the ${\rm ADO}_3$ invariant and becomes the original ${\rm ADO}_3$ polynomial after setting $t=-1$ and rescaling $x \mapsto \zeta_{3}^2 x$,
\begin{gather*} \zeta_{3} x^8 + \zeta_{3} x^7 + \big( \zeta_{3} - \zeta_{3}^{-1}\big) x^6 - \zeta_{3}^{-1}x^5 - \zeta_{3}^{-1}x^4 + \zeta_{3} x^2 + \zeta_{3} x + \big( \zeta_{3} - \zeta_{3}^{-1} \big) + (x \rightarrow 1/x ).
 \end{gather*}

\appendix

\section{Further examples}\label{appendixA}
We record ${\rm ADO}_4$ polynomials of torus knots obtained from the algorithm together with the results in Section~\ref{section3.3}:
\begin{gather*}
{\rm ADO}_4[T(2,23)]= -{\rm i} x^{33}-{\rm i} x^{32}-(1+{\rm i}) x^{31}-(1+{\rm i}) x^{30}-x^{29}-x^{28}-{\rm i} x^{25}-{\rm i} x^{24}\\
\hphantom{{\rm ADO}_4[T(2,23)]=}{}
-(1+{\rm i})x^{23} -(1+{\rm i}) x^{22} -x^{21}-x^{20}-{\rm i} x^{17}-{\rm i} x^{16}-(1+{\rm i}) x^{15}\\
\hphantom{{\rm ADO}_4[T(2,23)]=}{}-(1+{\rm i})x^{14}-x^{13}-x^{12} -{\rm i} x^{10}-2 {\rm i} x^9+(1-2 {\rm i}) x^8 -2 {\rm i} x^7 -{\rm i} x^6-x^4\\
\hphantom{{\rm ADO}_4[T(2,23)]=}{} -{\rm i} x^2-2 {\rm i}x+(1-2 {\rm i}) + (x \rightarrow 1/x ),
\\
{\rm ADO}_4[T(2,25)]= x^{36}+x^{35}+(1-{\rm i}) x^{34}+(1-{\rm i}) x^{33}-{\rm i} x^{32}-{\rm i} x^{31}+x^{28}+x^{27}\\
 \hphantom{{\rm ADO}_4[T(2,25)]=}{}
 +(1-{\rm i})x^{26}+(1-{\rm i}) x^{25}-{\rm i} x^{24}-{\rm i} x^{23}+x^{20}+x^{19}+(1-{\rm i}) x^{18}\\
 \hphantom{{\rm ADO}_4[T(2,25)]=}{}
 +(1-{\rm i}) x^{17}-{\rm i}x^{16}-{\rm i} x^{15} +x^{12}-{\rm i} x^{10} -2 {\rm i} x^9-(1+2 {\rm i}) x^8-2 {\rm i} x^7-{\rm i} x^6\\
 \hphantom{{\rm ADO}_4[T(2,25)]=}{}
 +x^4-{\rm i} x^2-2 {\rm i}x-(1+2 {\rm i}) + (x \rightarrow 1/x ),
\\
{\rm ADO}_4[T(2,27)]= {\rm i} x^{39}+ {\rm i} x^{38} + (1+{\rm i})x^{37} + (1+{\rm i})x^{36} + x^{35} + x^{34} + {\rm i} x^{31}+ {\rm i} x^{30}\\
\hphantom{{\rm ADO}_4[T(2,27)]=}{} + (1+{\rm i})x^{29} + (1+{\rm i})x^{28} + x^{27} + x^{26} + {\rm i} x^{23}+ {\rm i} x^{22} + (1+{\rm i})x^{21}\\
\hphantom{{\rm ADO}_4[T(2,27)]=}{} + (1+{\rm i})x^{20} + x^{19} + x^{18} + {\rm i} x^{15}+ {\rm i} x^{14} + (1+{\rm i})x^{13} + (1+{\rm i}2)x^{12} \\
\hphantom{{\rm ADO}_4[T(2,27)]=}{} + (1+{\rm i}) x^{11} + {\rm i} x^{10} + (-1+{\rm i})x^{9} - x^{8} + (-1+{\rm i})x^{7} + {\rm i} x^{6} + (1+{\rm i}) x^{5} \\
\hphantom{{\rm ADO}_4[T(2,27)]=}{} + (1+{\rm i}2) x^{4} + (1+{\rm i})x^{3} + {\rm i} x^{2} + (-1+{\rm i})x -1 + (x \rightarrow 1/x ),
\\
{\rm ADO}_4[T(2,29)]= - x^{42} - x^{41}+ (-1+{\rm i})x^{40} + (-1+{\rm i})x^{39} + {\rm i} x^{38} + {\rm i} x^{37} - x^{34} - x^{33}\\
\hphantom{{\rm ADO}_4[T(2,29)]=}{} + (-1+{\rm i})x^{32} + (-1+{\rm i})x^{31} + {\rm i} x^{30} + {\rm i} x^{29} - x^{26} - x^{25} + (-1+{\rm i})x^{24} \\
\hphantom{{\rm ADO}_4[T(2,29)]=}{} + (-1+{\rm i})x^{23} + {\rm i} x^{22} + {\rm i} x^{21} - x^{18} - x^{17} + (-1+{\rm i}) x^{16} + (-1+{\rm i})x^{15} \\
\hphantom{{\rm ADO}_4[T(2,29)]=}{} + {\rm i} x^{14} + (1+{\rm i}) x^{13} + x^{12} + (1+{\rm i})x^{11} + {\rm i}x^{10} + (-1+{\rm i})x^{9} + (-1+{\rm i}2)x^{8}\\
\hphantom{{\rm ADO}_4[T(2,29)]=}{} + (-1+{\rm i})x^{7} + {\rm i} x^{6} + (1+{\rm i})x^{5} + x^{4} + (1+{\rm i})x^{3} + {\rm i} x^{2} + (-1+{\rm i})x \\
\hphantom{{\rm ADO}_4[T(2,29)]=}{}
-1+{\rm i}2 + (x \rightarrow 1/x ),
\\
{\rm ADO}_4[T(2,31)]= -{\rm i} x^{45}-{\rm i} x^{44}-(1+{\rm i}) x^{43}-(1+{\rm i}) x^{42}-x^{41}-x^{40}-{\rm i} x^{37}-{\rm i} x^{36}\\
\hphantom{{\rm ADO}_4[T(2,31)]=}{} -(1+{\rm i})x^{35} -(1+{\rm i}) x^{34} -x^{33}-x^{32}-{\rm i} x^{29}-{\rm i} x^{28}-(1+{\rm i}) x^{27}\\
\hphantom{{\rm ADO}_4[T(2,31)]=}{} -(1+{\rm i})x^{26}-x^{25}-x^{24} -{\rm i} x^{21} -{\rm i} x^{20} -(1+{\rm i}) x^{19} -(1+{\rm i}) x^{18} -x^{17}\\
\hphantom{{\rm ADO}_4[T(2,31)]=}{} -x^{16}-{\rm i}x^{14}-2 {\rm i} x^{13} +(1-2 {\rm i}) x^{12} -2 {\rm i} x^{11}-{\rm i} x^{10}-x^8 -{\rm i} x^6-2 {\rm i} x^5 \\
\hphantom{{\rm ADO}_4[T(2,31)]=}{} +(1-2 {\rm i})x^4-2 {\rm i} x^3-{\rm i} x^2 -1 + (x \rightarrow 1/x ),
\\
{\rm ADO}_4[T(2,33)]= x^{48} + x^{47}+ (1-{\rm i})x^{46} + (1-{\rm i})x^{45}- {\rm i} x^{44}- {\rm i} x^{43} +x^{40} + x^{39}\\
\hphantom{{\rm ADO}_4[T(2,33)]=}{} + (1-{\rm i})x^{38} + (1-{\rm i})x^{37} - {\rm i} x^{36}- {\rm i} x^{35} + x^{32} + x^{31} + (1-{\rm i})x^{30} \\
\hphantom{{\rm ADO}_4[T(2,33)]=}{} + (1-{\rm i})x^{29} - {\rm i} x^{28}- {\rm i} x^{27} + x^{24} + x^{23} + (1-{\rm i})x^{22} + (1-{\rm i})x^{21}- {\rm i} x^{20}\\
\hphantom{{\rm ADO}_4[T(2,33)]=}{} - {\rm i} x^{19} + x^{16} - {\rm i} x^{14} -{\rm i}2 x^{13} + (-1-{\rm i}2 )x^{12}- {\rm i}2 x^{11} - {\rm i} x^{10} + x^{8} -{\rm i}x^{6}\\
\hphantom{{\rm ADO}_4[T(2,33)]=}{} -{\rm i}2x^{5} + (-1-{\rm i}2)x^{4} -{\rm i}2x^{3} -{\rm i}x^{2} + 1 + (x \rightarrow 1/x ),
\\
{\rm ADO}_4[T(2,35)]= {\rm i} x^{51}+{\rm i} x^{50}+(1+{\rm i}) x^{49}+(1+{\rm i}) x^{48}+x^{47}+x^{46}+{\rm i} x^{43}+{\rm i} x^{42}\\
\hphantom{{\rm ADO}_4[T(2,35)]=}{} +(1+{\rm i})x^{41} +(1+{\rm i}) x^{40}+x^{39} +x^{38}+{\rm i} x^{35}+{\rm i} x^{34}+(1+{\rm i}) x^{33}\\
\hphantom{{\rm ADO}_4[T(2,35)]=}{} +(1+{\rm i})x^{32}+x^{31}+x^{30} +{\rm i} x^{27} +{\rm i} x^{26} +(1+{\rm i}) x^{25}+(1+{\rm i}) x^{24} +x^{23}\\
\hphantom{{\rm ADO}_4[T(2,35)]=}{} +x^{22} +{\rm i}x^{19}+{\rm i} x^{18} +(1+{\rm i}) x^{17} +(1+2 {\rm i}) x^{16}+(1+{\rm i}) x^{15} +{\rm i} x^{14}\\
\hphantom{{\rm ADO}_4[T(2,35)]=}{} -(1-{\rm i})x^{13}-x^{12}-(1-{\rm i}) x^{11} +{\rm i} x^{10} +(1+{\rm i}) x^9 +(1+2 {\rm i}) x^8\\
\hphantom{{\rm ADO}_4[T(2,35)]=}{} +(1+{\rm i}) x^7+{\rm i} x^6-(1-{\rm i})x^5 -x^4-(1-{\rm i}) x^3 +{\rm i} x^2 +(1+{\rm i}) x \\
\hphantom{{\rm ADO}_4[T(2,35)]=}{} +(1+2 {\rm i}) + (x \rightarrow 1/x ),
\\
{\rm ADO}_4[T(2,37)]= -x^{54}-x^{53}-(1-{\rm i}) x^{52}-(1-{\rm i}) x^{51}+{\rm i} x^{50}+{\rm i} x^{49}-x^{46}-x^{45}\\
\hphantom{{\rm ADO}_4[T(2,37)]=}{} -(1-{\rm i})x^{44} -(1-{\rm i}) x^{43}+{\rm i} x^{42} +{\rm i} x^{41}-x^{38}-x^{37}-(1-{\rm i}) x^{36}\\
\hphantom{{\rm ADO}_4[T(2,37)]=}{} -(1-{\rm i}) x^{35}+{\rm i}x^{34}+{\rm i} x^{33} -x^{30} -x^{29} -(1-{\rm i}) x^{28}-(1-{\rm i}) x^{27} +{\rm i} x^{26}\\
\hphantom{{\rm ADO}_4[T(2,37)]=}{} +{\rm i}x^{25}-x^{22}-x^{21}-(1-{\rm i}) x^{20} -(1-{\rm i}) x^{19} +{\rm i} x^{18} +(1+{\rm i}) x^{17}+x^{16}\\
\hphantom{{\rm ADO}_4[T(2,37)]=}{} +(1+{\rm i})x^{15}+{\rm i} x^{14} -(1-{\rm i}) x^{13} -(1-2 {\rm i}) x^{12}-(1-{\rm i}) x^{11} +{\rm i} x^{10}\\
\hphantom{{\rm ADO}_4[T(2,37)]=}{} +(1+{\rm i})x^9+x^8+(1+{\rm i}) x^7+{\rm i} x^6-(1-{\rm i}) x^5 -(1-2 {\rm i}) x^4 -(1-{\rm i}) x^3\\
\hphantom{{\rm ADO}_4[T(2,37)]=}{} +{\rm i} x^2 +(1+{\rm i}) x+1 + (x \rightarrow 1/x ),
\\
{\rm ADO}_4[T(2,39)]= -{\rm i} x^{57}-{\rm i} x^{56}-(1+{\rm i}) x^{55}-(1+{\rm i}) x^{54}-x^{53}-x^{52}-{\rm i} x^{49}-{\rm i} x^{48}\\
\hphantom{{\rm ADO}_4[T(2,39)]=}{} -(1+{\rm i})x^{47} -(1+{\rm i}) x^{46} -x^{45}-x^{44}-{\rm i} x^{41}-{\rm i} x^{40}-(1+{\rm i}) x^{39}\\
\hphantom{{\rm ADO}_4[T(2,39)]=}{} -(1+{\rm i})x^{38}-x^{37}-x^{36} -{\rm i} x^{33} -{\rm i} x^{32} -(1+{\rm i}) x^{31} -(1+{\rm i}) x^{30} -x^{29}\\
\hphantom{{\rm ADO}_4[T(2,39)]=}{} -x^{28}-{\rm i}x^{25}-{\rm i} x^{24} -(1+{\rm i}) x^{23}-(1+{\rm i}) x^{22} -x^{21}-x^{20} -{\rm i} x^{18}-2 {\rm i} x^{17}\\
\hphantom{{\rm ADO}_4[T(2,39)]=}{} +(1-2 {\rm i})x^{16} -2 {\rm i} x^{15} -{\rm i} x^{14}-x^{12}-{\rm i} x^{10} -2 {\rm i} x^9+(1-2 {\rm i}) x^8 -2 {\rm i} x^7\\
\hphantom{{\rm ADO}_4[T(2,39)]=}{} -{\rm i} x^6-x^4 -{\rm i}x^2-2 {\rm i} x +(1-2 {\rm i}) + (x \rightarrow 1/x ).
\end{gather*}

\section[Comparison with the R-matrix approach]{Comparison with the $\boldsymbol{R}$-matrix approach}\label{appendixB}
We perform an independent computation of the ADO polynomial using its $R$-matrix formulation~\cite{ADO, BDGG} to strengthen the Conjecture~\ref{conjecture1}. We summarize the ingredients for the computation~\cite{BDGG}. A $(1,1)$-tangle diagram of $T(2,2s+1)$ consists of three kinds of building blocks: oriented caps, cups, and crossings, respectively,
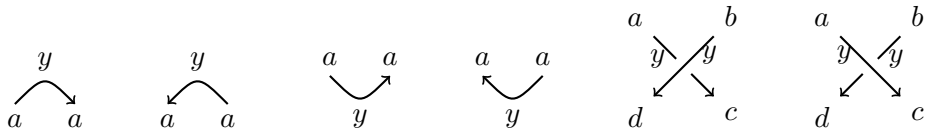
\begin{figure}[h!]\centering
\begin{tikzpicture}[scale=1]
	\draw[thick,->] (0,0) node[below]{$a$} .. controls (0.4,0.4) .. (0.8,0) node[below] {$a$}
	node[pos=0.5, above]{$y$};
	\end{tikzpicture}\qquad \begin{tikzpicture}
\draw[thick,->] (0.8,0) node[below]{$a$} .. controls (0.4,0.4) .. (0,0) node[below] {$a$}
node[pos=0.5, above]{$y$};
\end{tikzpicture} \qquad \begin{tikzpicture}
\draw[thick,->] (0,0) node[above]{$a$} .. controls (0.4,-0.4) .. (0.8,0) node[above] {$a$}
node[pos=0.5, below]{$y$};
\end{tikzpicture}\qquad \begin{tikzpicture}
\draw[thick,->] (0.8,0) node[above]{$a$} .. controls (0.4,-0.4) .. (0,0) node[above] {$a$}
node[pos=0.5, below]{$y$};
\end{tikzpicture}\qquad \begin{tikzpicture}
\draw[thick,->] (0.8,0) node[above right]{$b$} -- (0,-0.8) node[below left] {$d$};
\draw[thick] (0,0) node[above left]{$a$} -- (0.3,-0.3) node[left]{$y$};
\draw[thick,->] (0.5,-0.5) node[above right]{$y$} -- (0.8,-0.8) node[below right] {$c$};
\end{tikzpicture}\qquad \begin{tikzpicture}
\draw[thick,->] (0,0) node[above left]{$a$} -- (0.8,-0.8) node[below right] {$c$};
\draw[thick] (0.8,0) node[above right]{$b$} -- (0.5,-0.3) node[right]{$y$};
\draw[thick,->] (0.3,-0.5) node[above left]{$y$} -- (0,-0.8) node[below left] {$d$};
\end{tikzpicture}
\end{figure}
\begin{gather*}
\epsilon_{a}[y]=1,\qquad \epsilon^{\ast}_{a}[y]=q^{2a(r-1)}y^{1-r}, \qquad \eta_{a}[y]=1,\qquad \eta^{\ast}_{a}[y]=q^{2a(1-r)}y^{r-1},
\\
R^{a,b}_{c,d}[y]= \delta_{a-c, d-b} \theta_{a\geq c} \theta_{d\geq b} (-y)^{a-c} q^{(c-a)(a+b+1)+2cd} z y^{-d-c}\\
\hphantom{R^{a,b}_{c,d}[y]=}{}\times \frac{\big(q^{2(a-1)}/y^2 ; q^{-2}\big)_{a-c} \big(q^{2(b+1)}; q^{2}\big)_{a-c}}{\big(q^{-2} ; q^{-2}\big)_{a-c}},
\\
\big( R^{-1} \big)^{a,b}_{c,d}[y]= \delta_{a-c, d-b} \theta_{a\geq c} \theta_{d\geq b} (-y)^{a-c}q^{(c-a)(a+b+1)-2ab} z y^{b+a}\\
\hphantom{\big( R^{-1} \big)^{a,b}_{c,d}[y]=}{}\times \frac{\big(q^{2(a-1)}/y^2 ; q^{-2}\big)_{a-c} \big(q^{2(b+1)}; q^{2}\big)_{a-c}}{\big(q^{2} ; q^{2}\big)_{a-c}},
\\
\delta_{a,b} = \begin{cases}
 1, & a=b, \\
 0, & \text{otherwise},
 \end{cases}\qquad \theta_{a\geq b}= \begin{cases}
 1, & a\geq b, \\
 0, & \text{otherwise},
 \end{cases}
\end{gather*}
where $a$, $b$, $c$, $d$ are subset of variables $a_{1},\dots,a_{m}$ in the tangle diagram and $(w;q)_t$ is the $q$-Pochhammer symbol (see Section~\ref{section3.4}). The above formulas are in the same order as the diagrams. From these ingredients, a~function that gives rise to the ADO polynomial can be defined~as
\begin{gather*}
G^{\times}_{D}(q,y,z,r; a_{1},\dots,a_{m}) := d[y] \delta_{a_{1},0} \delta_{a_{m},0} \prod_{\rm crossings} R \prod_{\rm crossings} R^{-1} \prod_{\rm caps} \epsilon \prod_{\rm caps} \epsilon^{\ast} \prod_{\rm cups} \eta \prod_{\rm cups} \eta^{\ast},
\\
d[y]=\prod_{j=2}^{r} \frac{1}{q^{j}y - q^{-j}y^{-1}}=(-y)^{r-1}q^{\frac{1}{2}r(r+1)-1}\frac{1}{\big(q^4 y^2;q^2\big)_{r-1}}.
\end{gather*}
At $q=\zeta_{2r}$, $y=\zeta_{2r}^{\alpha}$, $z=\zeta_{2r}^{\alpha^2}$, an (unnormalized) ADO polynomial $N^{r}_K(\alpha)$ is
\begin{gather*}
N^{r}_K(\alpha) = G^{\times}_{D}\big(\zeta_{2r},\zeta_{2r}^{\alpha},\zeta_{2r}^{\alpha^2},r; a_{1},\dots,a_{m}\big).
\end{gather*}
The quantity computed in \cite{BDGG} is a normalized version
\begin{gather*}
\hat{N}^{r}_K (\alpha) := i^{1-r} \big(y^r - y^{-r}\big)N^{r}_K(\alpha -1).
\end{gather*}
The change of normalization between $\hat{N}^{r}_K (\alpha)$ and our ${\rm ADO}_p$ for zero framed knots is
\begin{gather*}
{\rm ADO}_p(x) \cong \frac{\hat{N}^{r=p}_K (\alpha ;y)}{y - y^{-1}}\bigg\rvert_{y\rightarrow x^{1/2},\, x\rightarrow cx} \cong \operatorname{num}[N^{r}_K(\alpha -1;y)]\bigg\rvert_{y\rightarrow x^{1/2},\, x\rightarrow cx},\qquad c \in \mathbb{C}^{\ast},
\end{gather*}
where $\cong$ denotes equivalence up to an overall monomial and an overall constant. The r.h.s.\ is due to the structure of $G^{\times}_{D}(r)$ such that $N^{r}_K(\alpha -1)$ always contains $\big(y - y^{-1}\big)/\big(y^r - y^{-r}\big)$ for any knot (for details see \cite[Section~2.4]{BDGG}). We denote the numerator of $N^{r}_K(\alpha -1;y)$ as $\operatorname{num}\big[N^{r}_K(\alpha -1;y)\big]$.

A $(1,1)$-tangle diagram of $T(2,2s+1)$, which consists of $(2s+1)$-crossings is
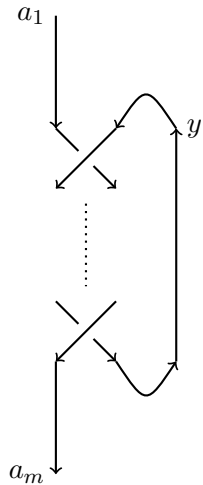
\begin{center}
\begin{tikzpicture}[scale=1]
\draw[thick,->] (0,0) node [left] {$a_1$} -- (0,-1.5);
\draw[thick] (0,-1.5) -- (0.3,-1.8);
\draw[thick,->] (0.5,-2) -- (0.8,-2.3);
\draw[thick,->] (0.8,-1.5) -- (0,-2.3);
\draw[thick] (0,-3.8) -- (0.3,-4.1);
\draw[thick,->] (0.5,-4.3) -- (0.8,-4.6);
\draw[thick,->] (0.8,-3.8) -- (0,-4.6);
\draw[thick,->] (0,-4.6) -- (0,-6.1) node [left] {$a_m$} ;
\draw[thick,->] (0.8,-4.6) .. controls (1.2,-5.2) .. (1.6,-4.6);
\draw[thick,->] (1.6,-4.6) -- (1.6,-1.5) node [right] {$y$};
\draw[thick,->] (1.6,-1.5) .. controls (1.2,-0.9) .. (0.8,-1.5);
\draw[dotted, thick] (0.4,-2.5) -- (0.4,-3.6);
\end{tikzpicture}
\end{center}
The vertical dots represent the same type of crossings. Applying the formula to the diagram, we have schematically
\begin{gather*}
G^{\times}_{D}(q,y,z,r; a_{1},\dots,a_{m}) = d[y] \delta_{a_{1},0} \delta_{a_{m},0} \left( \prod_{i=1}^{2s+1} R_{i} \right) \eta \epsilon^{\ast},
\end{gather*}
where $m=m(s)$. The ADO polynomials for $T(2,3)$ are listed in \cite[Appendix B]{BDGG}. Using the above relation ($c=1$), we find an agreement that
\begin{gather*}
\hat{N}^{3}_{T(2,3)}(\alpha ;y)= q^2 \big(y^5 - y^{-5}\big) + q\big(y - y^{-1}\big) \Rightarrow \frac{\hat{N}^{3}_K (\alpha ;y)}{y - y^{-1}}\bigg\rvert_{y\rightarrow x^{1/2}} \cong {\rm ADO}_3[T(2,3)](x),
\\
\hat{N}^{4}_{T(2,3)}(\alpha ;y)= q^2 \big(y^7 - y^{-7}\big) + \big(y^3 - y^{-3}\big) + q^2\big(y - y^{-1}\big) \Rightarrow \frac{\hat{N}^{4}_K (\alpha ;y)}{y - y^{-1}}\bigg\rvert_{y\rightarrow x^{1/2}} \\
\hphantom{\hat{N}^{4}_{T(2,3)}(\alpha ;y)}{}
\cong {\rm ADO}_4[T(2,3)](x).
\end{gather*}
We next check $T(2,5)$ case. The computation of $G^{\times}_{D}$ yields
\begin{gather*}
\operatorname{num}\big[N^{3}_{T(2,5)}(\alpha -1)\big] = -\sqrt[3]{-1} y^8 -\sqrt[3]{-1}y^6 -\frac{1}{2} \sqrt[3]{-1} \big(3-{\rm i} \sqrt{3}\big)y^4 -\frac{1}{2} \sqrt[3]{-1}\big(1-{\rm i} \sqrt{3}\big) y^2 \\
\hphantom{\operatorname{num}\big[N^{3}_{T(2,5)}(\alpha -1)\big] =}{}
 -\frac{1}{2} \sqrt[3]{-1} \big(1-{\rm i} \sqrt{3}\big) + \left( y \rightarrow \frac{1}{y} \right) \bigg\rvert_{y\rightarrow x^{1/2}} \cong {\rm ADO}_3[T(2,5)](x).
\end{gather*}
We now list several more verifications of the ${\rm ADO}_3$ formula in Section~\ref{section3.1}:
\begin{gather*}
\operatorname{num}\big[N^{3}_{T(2,7)}(\alpha -1)\big] = -\sqrt[3]{-1} y^{12} -\sqrt[3]{-1}y^{10}-\frac{1}{2} \sqrt[3]{-1} \big(3-{\rm i} \sqrt{3}\big) y^8 -\frac{1}{2}\sqrt[3]{-1} \big(1-{\rm i} \sqrt{3}\big) y^6\\
\hphantom{\operatorname{num}\big[N^{3}_{T(2,7)}(\alpha -1)\big] =}{} -\frac{1}{2} \sqrt[3]{-1} \big(1-{\rm i} \sqrt{3}\big)y^4 -\sqrt[3]{-1} + \left( y \rightarrow \frac{1}{y} \right) \bigg\rvert_{y\rightarrow x^{1/2}} \\
\hphantom{\operatorname{num}\big[N^{3}_{T(2,7)}(\alpha -1)\big]}{}
\cong {\rm ADO}_3[T(2,7)](x),
\\
\operatorname{num}\big[N^{3}_{T(2,9)}(\alpha -1)\big] = -\frac{1}{2} \sqrt[6]{-1} \big(\sqrt{3}+{\rm i}\big) y^{16} -\frac{1}{2} \sqrt[6]{-1}\big(\sqrt{3}+{\rm i}\big)y^{14} -\sqrt[6]{-1} \sqrt{3}y^{12}\\
\hphantom{\operatorname{num}\big[N^{3}_{T(2,9)}(\alpha -1)\big] =}{}
 -\frac{1}{2} \sqrt[6]{-1}\big(\sqrt{3}-{\rm i}\big) y^{10} -\frac{1}{2} \sqrt[6]{-1} \big(\sqrt{3}-{\rm i}\big) y^8 -\frac{1}{2} \sqrt[6]{-1} \big(\sqrt{3}+{\rm i}\big)y^4\\
\hphantom{\operatorname{num}\big[N^{3}_{T(2,9)}(\alpha -1)\big] =}{}
 -\frac{1}{2} \sqrt[6]{-1}\big(\sqrt{3}+{\rm i}\big) y^2 -\sqrt[6]{-1} \sqrt{3} + \left( y \rightarrow \frac{1}{y} \right) \bigg\rvert_{y\rightarrow x^{1/2}}\\
\hphantom{\operatorname{num}\big[N^{3}_{T(2,9)}(\alpha -1)\big]}{}
 \cong {\rm ADO}_3[T(2,9)](x),
\\
\operatorname{num}\big[N^{3}_{T(2,11)}(\alpha -1)\big] = -(-1)^{2/3} y^{20} -(-1)^{2/3}y^{18} -\frac{1}{2} (-1)^{2/3}\big(3-{\rm i} \sqrt{3}\big) y^{16}\\
\hphantom{\operatorname{num}\big[N^{3}_{T(2,11)}(\alpha -1)\big] =}{} - \frac{1}{2} (-1)^{2/3} \big(1-{\rm i} \sqrt{3}\big) y^{14} -\frac{1}{2} (-1)^{2/3} \big(1-{\rm i} \sqrt{3}\big) y^{12}\\
 \hphantom{\operatorname{num}\big[N^{3}_{T(2,11)}(\alpha -1)\big] =}{}
 -(-1)^{2/3}y^8 -(-1)^{2/3} y^6 -\frac{1}{2} (-1)^{2/3} \big(3-{\rm i} \sqrt{3}\big) y^4\\
\hphantom{\operatorname{num}\big[N^{3}_{T(2,11)}(\alpha -1)\big] =}{}
 -\frac{1}{2} (-1)^{2/3} \big(1-{\rm i} \sqrt{3}\big)y^2 -\frac{1}{2} (-1)^{2/3}\big(1-{\rm i} \sqrt{3}\big)+ \left( y \rightarrow \frac{1}{y} \right)\! \bigg\rvert_{y\rightarrow x^{1/2}}\\
\hphantom{\operatorname{num}\big[N^{3}_{T(2,11)}(\alpha -1)\big]}{} \cong {\rm ADO}_3[T(2,11)](x),
\\
\operatorname{num}\big[N^{3}_{T(2,13)}(\alpha -1)\big] = \frac{1}{2} \big(-1-{\rm i} \sqrt{3}\big) y^{24} +\frac{1}{2} \big({-}1-{\rm i} \sqrt{3}\big) y^{22} +\frac{1}{2} \big({-}3-i \sqrt{3}\big) y^{20} -y^{18}\\
\hphantom{\operatorname{num}\big[N^{3}_{T(2,13)}(\alpha -1)\big] =}{}
 -y^{16} +\frac{1}{2} \big({-}1-{\rm i} \sqrt{3}\big) y^{12} +\frac{1}{2} \big({-}1-{\rm i} \sqrt{3}\big) y^{10} +\frac{1}{2} \big({-}3-{\rm i} \sqrt{3}\big)y^8\\
\hphantom{\operatorname{num}\big[N^{3}_{T(2,13)}(\alpha -1)\big] =}{}
-y^6 -y^4 +\frac{1}{2} \big({-}1-{\rm i}\sqrt{3}\big)+ \left( y \rightarrow \frac{1}{y} \right) \bigg\rvert_{y\rightarrow x^{1/2}}\\
\hphantom{\operatorname{num}\big[N^{3}_{T(2,13)}(\alpha -1)\big]}{}
\cong {\rm ADO}_3[T(2,13)](x),
\\
\operatorname{num}\big[N^{3}_{T(2,15)}(\alpha -1)\big] = \frac{1}{2} \big(1-{\rm i} \sqrt{3}\big) y^{28} +\frac{1}{2} \big(1-{\rm i} \sqrt{3}\big) y^{26} -{\rm i} \sqrt{3}y^{24} +\frac{1}{2} \big({-}1-{\rm i} \sqrt{3}\big)y^{22}\\
\hphantom{\operatorname{num}\big[N^{3}_{T(2,15)}(\alpha -1)\big] =}{}
 +\frac{1}{2} \big({-}1-{\rm i} \sqrt{3}\big)y^{20} +\frac{1}{2} \big(1-{\rm i} \sqrt{3}\big) y^{16} +\frac{1}{2} \big(1-{\rm i} \sqrt{3}\big)y^{14} -i \sqrt{3} y^{12}\\
\hphantom{\operatorname{num}\big[N^{3}_{T(2,15)}(\alpha -1)\big] =}{}
 +\frac{1}{2} \big({-}1-{\rm i} \sqrt{3}\big) y^{10} +\frac{1}{2} \big({-}1-{\rm i} \sqrt{3}\big) y^8 +\frac{1}{2} \big(1-{\rm i} \sqrt{3}\big) y^4\\
\hphantom{\operatorname{num}\big[N^{3}_{T(2,15)}(\alpha -1)\big] =}{}
 +\frac{1}{2} \big(1-{\rm i} \sqrt{3}\big) y^2 -{\rm i} \sqrt{3} + \left( y \rightarrow \frac{1}{y} \right) \bigg\rvert_{y\rightarrow x^{1/2}}
\cong {\rm ADO}_3[T(2,15)](x),
\\
\operatorname{num}\big[N^{3}_{T(2,17)}(\alpha -1)\big] = \frac{1}{2} {\rm i} \big(\sqrt{3}+{\rm i}\big) y^{32} +\frac{1}{2}{\rm i} \big(\sqrt{3}+{\rm i}\big) y^{30} + {\rm i}\sqrt{3} y^{28} +\frac{1}{2} \big(1+{\rm i} \sqrt{3}\big) y^{26}\\
\hphantom{\operatorname{num}\big[N^{3}_{T(2,17)}(\alpha -1)\big] =}{}
 +\frac{1}{2} \big(1+{\rm i} \sqrt{3}\big) y^{24} +\frac{1}{2} {\rm i} \big(\sqrt{3}+{\rm i}\big) y^{20} +\frac{1}{2} {\rm i}\big(\sqrt{3}+{\rm i}\big) y^{18} +{\rm i} \sqrt{3}y^{16}\\
\hphantom{\operatorname{num}\big[N^{3}_{T(2,17)}(\alpha -1)\big] =}{}
+\frac{1}{2} \big(1+{\rm i} \sqrt{3}\big)y^{14} +\frac{1}{2} \big(1+{\rm i} \sqrt{3}\big) y^{12} +\frac{1}{2} {\rm i} \big(\sqrt{3}+{\rm i}\big)y^8\\
\hphantom{\operatorname{num}\big[N^{3}_{T(2,17)}(\alpha -1)\big] =}{}
+\frac{1}{2} {\rm i} \big(\sqrt{3}+{\rm i}\big)y^6+{\rm i} \sqrt{3} y^4+\frac{{\rm i}\sqrt{3}}{y^4} +\frac{1}{2} \big(1+{\rm i} \sqrt{3}\big) y^2\\
\hphantom{\operatorname{num}\big[N^{3}_{T(2,17)}(\alpha -1)\big] =}{}
+\frac{1}{2} \big(1+{\rm i} \sqrt{3}\big)+ \left( y \rightarrow \frac{1}{y} \right) \bigg\rvert_{y\rightarrow x^{1/2}} \cong {\rm ADO}_3[T(2,17)](x).
\end{gather*}

We next verify ${\rm ADO}_4$ polynomials:
\begin{gather*}
\operatorname{num}\big[N^{4}_{T(2,7)}(\alpha -1)\big] = -\sqrt[3]{-1} y^{12}-\sqrt[3]{-1}y^{10}-\frac{1}{2} \sqrt[3]{-1} \big(3-{\rm i}\sqrt{3}\big) y^8-\frac{1}{2} \sqrt[3]{-1} \big(1-{\rm i} \sqrt{3}\big) y^6\\
\hphantom{\operatorname{num}\big[N^{4}_{T(2,7)}(\alpha -1)\big] =}{}
 -\frac{1}{2} \sqrt[3]{-1}\big(1-{\rm i} \sqrt{3}\big) y^4-\sqrt[3]{-1}+ \left( y \rightarrow \frac{1}{y} \right) \bigg\rvert_{y\rightarrow x^{1/2}}\\
\hphantom{\operatorname{num}\big[N^{4}_{T(2,7)}(\alpha -1)\big]}{}
 \cong {\rm ADO}_4[T(2,7)](x),
\\
\operatorname{num}\big[N^{4}_{T(2,9)}(\alpha -1)\big] = -\sqrt[4]{-1} y^{24}-\sqrt[4]{-1}y^{22}-(1-{\rm i}) \sqrt[4]{-1} y^{20}-(1-{\rm i}) \sqrt[4]{-1} y^{18}\\
\hphantom{\operatorname{num}\big[N^{4}_{T(2,9)}(\alpha -1)\big] =}{}
+(-1)^{3/4} y^{16} +(-1)^{3/4}y^{14}-\sqrt[4]{-1} y^8 +(-1)^{3/4} y^4+2 (-1)^{3/4} y^2\\
\hphantom{\operatorname{num}\big[N^{4}_{T(2,9)}(\alpha -1)\big] =}{}
 +(1+2 {\rm i}) \sqrt[4]{-1} + \left( y \rightarrow \frac{1}{y} \right) \bigg\rvert_{y\rightarrow x^{1/2}}\cong {\rm ADO}_4[T(2,9)](x),
\\
\operatorname{num}\big[N^{4}_{T(2,11)}(\alpha -1)\big] = {\rm i} y^{30}+{\rm i} y^{28}+(1+{\rm i})y^{26}+(1+{\rm i})y^{24}+y^{22}+y^{20}+{\rm i}y^{14}+{\rm i} y^{12}\\
\hphantom{\operatorname{num}\big[N^{4}_{T(2,11)}(\alpha -1)\big] =}{}
+(1+{\rm i})y^{10} +(1+2 {\rm i}) y^8 +(1+{\rm i}) y^6+iy^4-(1-{\rm i}) y^2 -1\\
\hphantom{\operatorname{num}\big[N^{4}_{T(2,11)}(\alpha -1)\big] =}{}
 + \left( y \rightarrow \frac{1}{y} \right) \bigg\rvert_{y\rightarrow x^{1/2}} \cong {\rm ADO}_4[T(2,11)](x),
\\
\operatorname{num}\big[N^{4}_{T(2,13)}(\alpha -1)\big] = (-1+{\rm i}) y^{36}-(1-{\rm i}) y^{34}+2 {\rm i}y^{32}+2 {\rm i} y^{30}+(1+{\rm i})y^{28}+(1+{\rm i}) y^{26}\\
 \hphantom{\operatorname{num}\big[N^{4}_{T(2,13)}(\alpha -1)\big] =}{}
 -(1-{\rm i})y^{20} -(1-{\rm i}) y^{18}+2{\rm i}y^{16}+2 {\rm i} y^{14}+(1+{\rm i})y^{12}+2 y^{10}\\
 \hphantom{\operatorname{num}\big[N^{4}_{T(2,13)}(\alpha -1)\big] =}{}
 +(1-{\rm i})y^8 +2 y^6 +(1+{\rm i}) y^4 +2 {\rm i}y^2 +(1+3 {\rm i}) + \left( y \rightarrow \frac{1}{y} \right)\! \bigg\rvert_{y\rightarrow x^{1/2}}\!\\
 \hphantom{\operatorname{num}\big[N^{4}_{T(2,13)}(\alpha -1)\big]}{}
 \cong {\rm ADO}_4[T(2,13)](x).
\end{gather*}

\vspace{-4mm}

\subsection*{Acknowledgments}
I am grateful to Sergei Gukov for his valuable suggestions on a draft of this paper. I would like thank Angus Gruen for helpful conversations. I am also grateful to the referees for many helpful suggestions.

\pdfbookmark[1]{References}{ref}
\LastPageEnding

\end{document}